\documentclass{article}
\usepackage{spconf,amsmath,amsfonts,graphicx}
\usepackage[footnotesize]{caption}
\usepackage{url}

\newcommand{\bs}{\boldsymbol}
\newcommand{\Rbb}{\mathbb{R}}
\newcommand{\cl}{\mathcal}
\newcommand{\bb}{\mathbb}
\newcommand{\bbb}{|\!|\!|}
\newcommand{\ie}{\emph{i.e.}, }

\newcommand{\Id}{{\bf I}}
\DeclareMathOperator*{\argmin}{\arg\min}
\DeclareMathOperator{\prox}{prox}

\newcommand{\sq}{\vspace{-2mm}}

\title{\sq \sq Robust Phase Unwrapping by Convex Optimization\sq}
\name{Adriana Gonz{\'a}lez and Laurent Jacques
\thanks{Part of the work is funded by the DETROIT project (WIST3), convention no. 1017073, Walloon
  Region, Belgium. LJ is supported by the Belgian FRS-FNRS fund.}\sq\sq}
\address{\ninept ICTEAM, Universit\'e catholique de Louvain, Belgium}

\begin{document}

\maketitle

\begin{abstract}
The 2-D phase unwrapping problem aims at retrieving a ``phase'' image from its modulo $2\pi$ observations. Many applications, such as interferometry or synthetic aperture radar imaging, are concerned by this problem since they proceed by recording complex or modulated data from which a ``wrapped'' phase is extracted. Although 1-D phase unwrapping is trivial, a challenge remains in higher dimensions to overcome two common problems: noise and discontinuities in the true phase image. In contrast to state-of-the-art techniques, this work aims at simultaneously unwrap and denoise the phase image. We propose a robust convex optimization approach that enforces data fidelity constraints expressed in the corrupted phase derivative domain while promoting a sparse phase prior. The resulting optimization problem is solved by the Chambolle-Pock primal-dual scheme. We show that under different observation noise levels, our approach compares favorably to those that perform the unwrapping and denoising in two separate steps.
\end{abstract}

\begin{keywords}
\noindent \ninept Phase unwrapping, convex optimization, chambolle-pock algorithm, sparse prior.
\end{keywords}

\section{Introduction}
\label{sec:intro}

The information contained in the ``phase image'' is essential in many applications such as Magnetic Resonance Imaging (MRI) \cite{Tian2008MRI} and interferometric synthetic aperture radar (InSAR) \cite{Rosen2000}. In these techniques the phase is not observed directly but computed from a complex signal. Therefore, the measured values are wrapped in the interval $[-\pi,\pi)$ and the observed signal presents $2\pi$ jumps. Phase unwrapping is the procedure that allows us to practically remove these discontinuities to obtain the actual phase image. 

Although, the modulo operation is quite trivial, its inversion can be very hard to solve. Phase unwrapping techniques need to be able to overcome, among other problems, discontinuities, noise and under-sampling of the phase.

In one dimension, the unwrapping process is straightforward since there is only one possible ``path'' and a simple integration can recover the actual phase. However, this only works when the Itoh smoothness condition \cite{Ying2006_book} is satisfied, \ie when the absolute value of the phase gradient is lower or equal to $\pi$. The presence of noise or discontinuities could violate this condition, causing some unwrapping errors. 

Most existing methods extend this integration principle to two dimensions (2-D) and are denominated path-following algorithms. The problem in 2-D is the error propagation when the smoothness condition is violated. This occurs because the integration results depend on the chosen integration path and on the start and end points. The challenge remains in distinguishing jumps due to phase wrapping from those due to noise and discontinuities in the actual function. Several works have considered additional information such as pixel quality maps \cite{Arevalillo2002,Lu2007} to appropriately update the integration path. However, such algorithms have some difficulties to deal with high levels of noise.

In addition to these algorithms, several efforts have been made in the development of path-independent methods. State-of-the-art techniques rely on the global minimization of an energy function based on the classical $\ell_p$ norm of the error \cite{Ghiglia1998} or on a generalized $\ell_p$ norm \cite{Bioucas2007}. When minimizing the classical $\ell_p$ norm of the error as in \cite{Ghiglia1998}, the retrieved phase is generally smooth and sensitive to noise. In \cite{Bioucas2007} the authors prove that by minimizing a generalized $\ell_p$ norm using graph-cut techniques they are able to obtain an exact phase recovery (in noiseless scenarios). The developed algorithm is denoted PUMA and is considered a state-of-the-art method in phase unwrapping without denoising. Some works have included PUMA for unwrapping and they have added a denoising step either before \cite{Bioucas2008} or after \cite{Valadao2009} the unwrapping step. Both approaches have proved to outperform the classical $\ell_p$ norm minimization. While the approach of denoising before unwrapping \cite{Bioucas2008} deals better with high noise scenarios without discontinuities, it was reported that discontinuities are better preserved when the denoising is performed after the unwrapping \cite{Valadao2009}. 

We aim at developing a general numerical reconstruction method that is able to simultaneously unwrap and denoise the observed phase. We propose a convex optimization approach based on the minimization of a sparsity prior on the phase image under a data fidelity constraint expressed in the phase derivative domain and adapted to Gaussian distributed noise. Moreover, to stabilize the algorithm, the first phase component is assumed to be zero. The problem is solved by means of the primal-dual algorithm proposed by Chambolle and Pock \cite{ChambollePock_2011}. The results are compared with the state-of-the-art technique PUMA, whose output is denoised in a post-processing step. Although PUMA is faster and provides better results for the noiseless scenario, the proposed convex method provides better reconstruction quality for different noisy scenarios.

\section{Discrete forward model}
\label{sec:model}

Our work is concerned by the reconstruction of a phase image $\bs x \in \Rbb^N$ discretized over $N$ pixels on a regular Cartesian grid. The phase measurement process can be defined via the centralized modulo operator, or \emph{wrapping}, $\cl W$:
\begin{equation}
\label{eq:model}
	\bs y = \cl W(\bs x + \bs n),
\end{equation}
where $\cl W(\lambda) := [((\lambda + \pi) \mod 2\pi) - \pi] \in [-\pi,\pi) $ represents the component-wise wrapping, $\bs y \in \Rbb^{N}$ is the observed wrapped phase and $\bs n \in \Rbb^{N}$ stands for an additive phase sensing noise. 

Because of the modulo operation, we know that the actual phase $\bs x$ can be expressed as the sum of the observations $\bs y$ with integer multiples of $2\pi$, \ie $\bs x = \bs y + 2\pi\bs k$ for some $\bs k \in \bb Z^N$.

Interestingly, if the Itoh smoothness condition is satisfied and there is no noise in the measurement ($\bs n = \bs 0$), the image gradient can be indirectly observed through the following relation \cite{Ying2006_book, Valadao2009}:
\begin{equation}
\label{eq:grad_relation}
	\bs q := \cl W(\bs \nabla \bs y) = \bs \nabla \bs x,
\end{equation}
where $\bs \nabla : \Rbb^N \rightarrow \Rbb^{2N}$ is the (finite difference) gradient operator.

In the presence of low levels of noise such that the noisy phase satisfies the smoothness condition ($|(\bs \nabla (\bs x + \bs n))_j| \leq \pi$), \ie when $|(\bs \nabla \bs n)_j| \leq \pi - |(\bs \nabla\bs x)_j|$ on every pixel $j$, the relation in \eqref{eq:grad_relation} becomes: 
\begin{equation}
\label{eq:grad_relation2}
	\bs q := \cl W(\bs \nabla \bs y) = \bs \nabla (\bs x + \bs n).
\end{equation}
Assuming an additive white Gaussian noise model, \ie  $n_j \sim_{\rm iid} \cl{N}(0,\sigma^2)$, the Itoh condition is satisfied on all pixels if
$
\|\bs \nabla \bs n\|_{\infty} \leq \pi - \|\bs \nabla\bs x\|_{\infty}.
$
Since $\|\bs n\|_{\infty} = O(\sigma \sqrt{\log N})$ \cite{ledoux1991probability}, we have $\|\bs \nabla \bs n\|_{\infty} = O(\sigma \sqrt{\log N})$ and $
\|\bs \nabla\bs x\|_{\infty} \leq \pi - C \sigma \sqrt{\log N},
$
for some $C > 0$. This shows that the Itoh condition is increasingly harder to satisfy mainly when $\sigma$ increases.

Note that for any $\sigma$, the $\ell_2$ norm of $\bs n$ can be bounded using the Chernoff-Hoeffding bound \cite{hoeffding1963probability}: 
$$\| \bs n \|_2 \leq \varepsilon_{\rm n} := \sigma \sqrt{N + c\sqrt{N}},$$
which holds with high probability for $c = \cl O(1)$. 

When we encounter higher levels of noise and phase discontinuities, the relation in \eqref{eq:grad_relation2} is no longer valid. The error between the indirect observations and the phase gradient occurs in the $j$ pixels where the Itoh condition is not satisfied. Since $j$ represents a small percentage of the total amount of pixels (increasing with $\sigma$ and with the discontinuities), we can assume the error is sparse and it can be bounded using the $\ell_1$-norm:
\begin{equation}
\label{eq:fidelity}
	\| \bs q - \bs \nabla (\bs x + \bs n) \|_1 \leq \varepsilon_{\rm w},
\end{equation}
where $\varepsilon_{\rm w}$ is the error bound. In this paper, we estimate this error by computing $\| \bs q - \bs \nabla \bs x_{\rm p}\|_1$  using the reconstruction $\bs x_{\rm p}$ obtained by the PUMA method\footnote{Matlab code: \url{http://www.lx.it.pt/~bioucas/code.htm}} \cite{Bioucas2007}. This algorithm is known for estimating $\bs x + \bs n$ without denoising.

\section{Convex optimization approach}
\label{sec:prior}

It is rather clear that the ensemble of images having the same wrapped observations $\bs y$ forms a non-convex set: even when $\bs n = \bs 0$, if both $\bs x_1$ and $\bs x_2$ produce the same wrapped observation $\bs y$, there exist two integer vectors $\bs k_1$ and $\bs k_2$ such that $\bs y = \bs x_1 + 2\pi \bs k_1 = \bs x_2 + 2\pi \bs k_2$ and obviously $\bs x' = \lambda \bs x_1 + (1-\lambda) \bs x_2$ does not satisfy such a relation for all $\lambda \in [0,1]$ since $\lambda \bs k_1 + (1-\lambda) \bs k_2 \notin \bb Z^M$ for most $\lambda$. In front of such a discrete formalism, several authors have proposed methods based on combinatorial optimization and graph-cut techniques \cite{Valadao2009,Bioucas2007}. 

In this paper, we follow a different approach. We propose to relax the problem and to solve it using convex optimization by leveraging the differential relation \eqref{eq:grad_relation2}. Despite the wrapping operation, we expect that \eqref{eq:fidelity} holds for the phase signal $\bs x$ given an appropriate value $\varepsilon_{\rm w}$.  

Moreover, in order to circumvent the ill-conditioning of the problem in the presence of noise and ``non-Itoh'' phase discontinuities, we also regularized our method by an appropriate wavelet analysis prior model based on the structure of the phase, \ie a common tool used in many denoising methods \cite{donoho1995}. More specifically, we assume that the unwrapped phase image has a sparse or compressible representation in an orthonormal basis $\bs \Psi \in \Rbb^{N \times N}$, \ie the coefficients vector $\bs \Psi^T \bs x$ has few important values and its $\ell_1$-norm is expected to be small. The regularization process can then proceed by promoting a small $\ell_1$-norm in the wavelet projection of the phase image. The rationale of this is also to prevent fake phase jump reconstruction (since these increase locally the wavelet coefficient values) and to enforce the noise canceling. We also follow a common practice in the field which removes the (unsparse) scaling coefficients from the $\ell_1$-norm computation. 

Since both the differential fidelity and the wavelet prior are blind to the addition of a global constant, there is an ambiguity to estimate the phase up to such addition. In order to avoid this incertitude and to stabilize the convex optimization, we arbitrarily enforce the first phase component to be zero. We should note that the initial problem is itself ill-posed since, even if we solve it using \eqref{eq:model} directly with a perfect data prior model, a ``good'' solution would be determined up to a global addition of a multiple of $2\pi$. This constraint also induces the uniqueness of the solution.

%We also follow a common practice in the field which remove the (unsparse) scaling coefficients from the $\ell_1$-norm computation. However, in order to avoid the ambiguity to estimate the phase up to a global constant addition (and stabilizing the convex optimization), \ie both the differential fidelity and the wavelet prior being blind to any such addition, we arbitrarily enforce the first phase component to be zero. Remark that the initial problem is itself ill-posed since, even if we solved it from \eqref{eq:model} directly with a perfect data prior model, a ``good'' solution would be determined up to a global addition by a multiplie of $2\pi$. This constraint also induces the uniqueness of the solution.

Finally, since the noise level is assumed to be known and the noise $\ell_2$-norm to be bounded, we also propose to explicitly recover the noise part in an additive model where the unknown phase and noise are summed up to faithfully satisfy~\eqref{eq:fidelity}.

Gathering all these aspects, the proposed reconstruction program reads
\begin{align}
\argmin_{\bs u, \bs v \in \Rbb^N}\| \bs \Psi^T \bs u\|_1 &\ {\rm s.t.}\ \begin{cases}\| \bs v \|_2 \leq \varepsilon_{\rm n}\\ \| \bs q - \bs \nabla (\bs u + \bs v) \|_1 \leq \varepsilon_{\rm w}\\ u_1 = 0.
\label{eq:tvl1_1}
\end{cases}
\end{align}

Anticipating the scope of the next section, we may notice that by forming the vector $\bs w = (\bs u^T,\bs v^T)^T$, the convex minimization described above can be recast as\sq
\begin{multline}
\label{eq:tvl1_2}
\argmin_{\bs w \in \Rbb^{2N}} \|\bs \Psi^T \bs S_u \bs w\|_1 + \imath_{\cl C_1}(\bs S_v \bs w)\\[-2mm]
 + \imath_{\cl C_2}(\bs \nabla\,(\Id, \Id)\,\bs w) + \imath_{\Omega}(\bs w),
\end{multline}
where $\bs S_u, \bs S_v \in \bb R^{N\times 2N}$ are the selection operators of the first and the last $N$ elements of a vector in $\Rbb^{2N}$, respectively; $\Id \in \bb R^{N\times N}$ is the identity matrix; $\imath_{\cl P}(\bs x)$ is the (convex) indicator function of a convex set $\cl P$, which equals to $0$ if $\bs x\in \cl P$ and to $+\infty$ otherwise; and with the convex sets $\cl C_1=\{\bs z\in\Rbb^N:\| \bs z \|_2 \leq \varepsilon_{\rm n}\}$, $\cl C_2=\{\bs z\in\Rbb^{2N}:\| \bs q - \bs z \|_1 \leq \varepsilon_{\rm w}\}$ and $\Omega = \{ \bs z \in \Rbb^{2N}: z_1 = 0 \}$. In the next section we present the algorithm to solve \eqref{eq:tvl1_2} numerically.\sq

\section{Phase Unwrapping and Denoising Algorithm}
\label{sec:algorithm}

We are interested in finding the phase candidate that minimizes \eqref{eq:tvl1_2}, a problem that contains the sum of four lower semicontinuous convex functions from $\Rbb^D$ to $\Rbb \cup \{+\infty\}$, \ie they belong to the space $\Gamma_0(\Rbb^D)$ for some dimension $D\in \{N,2N\}$ \cite{combettes2011proximal}. In particular, we aim at solving the general optimization\sq
\begin{equation}
\label{eq:general-optim}
\min_{\bs w \in \bb R^{2N}} \sum_{j=1}^{p} F_j(\bs K_j \bs w) + H(\bs w),
\end{equation}
with $\bs K_j : \Rbb^N \rightarrow \Rbb^{W_j}$ and $p+1$ the number of
convex functions. For this, we use the Chambolle-Pock primal-dual algorithm defined in a Product Space \cite{ChambollePock_2011,IPI14}. As any other proximal algorithm \cite{parikh2013proximal}, this one relies on the definition of the proximal operator\sq 
\begin{equation*}
\prox_{\varphi} \bs z := \argmin_{\bs u\in\Rbb^D} \varphi(\bs u) + \tfrac{1}{2}\|\bs u - \bs z\|^2\sq
\end{equation*}
that is uniquely defined for any $\varphi \in \Gamma_0(\Rbb^D)$ for some $D\in \bb N$ \cite{combettes2011proximal}. By writing $\bs K={\rm diag}(\bs K_1, \cdots, \bs K_p)$, the CP iterations are
\begin{equation}
\label{eq:CP_Algorithm_PS}
\!\!\begin{cases}
\bs s_j^{(k+1)}\hspace{-2.5mm}&=\prox_{\nu F_j^\star} \big(\bs s_j^{{(k)}} + \nu \bs K_j \bar{\bs w}^{(k)} \big),\ j\in \{1,\cdots p\}\\
\bs w^{{(k+1)}}\hspace{-2.5mm}&= \prox_{\tfrac{\mu}{p} H} (\bs w^{{(k)}} - \tfrac{\mu}{p} \sum_{j=1}^p \bs K_j^* \bs s_j^{{(k+1)}}),\\
\bar{\bs w}^{(k+1)}\hspace{-2.5mm}&= 2\,\bs w^{(k+1)} - \bs w^{(k)},  
\end{cases}\sq
\end{equation}
with $\bs w^{k}$ tending to a minimizer $\bs w^*$ of \eqref{eq:general-optim} with $k\to + \infty$.  
To match the formulation \eqref{eq:general-optim} with the problem at hand \eqref{eq:tvl1_2}, we set $p=3$, $F_1 (\bs s_1) = \| \bs s_1\|_1$ and $F_2(\bs s_2) = \imath_{\cl C_1}(\bs s_2)$ for $\bs s_1, \bs s_2 \in \Rbb^N$; $F_3(\bs s_3) = \imath_{\cl C_2}(\bs s_3)$ for $\bs s_3 \in \Rbb^{2N}$; $H(\bs w) = \imath_{\Omega} (\bs w)$; $\bs K_1 = \bs \Psi^T \bs S_u$; $\bs K_2 = \bs S_v$ and $\bs K_3 = \bs \nabla\, (\Id, \Id)$. 

In order to apply the algorithm in \eqref{eq:CP_Algorithm_PS}, we must compute the proximal operators of $F_1^\star$, $F_2^\star$, $F_3^\star$ and $H$, the first three functions being the Legendre-Fenchel conjugate of their unstarred version. The proximal operator of $F^\star$ is determined via the one of $F$ thanks to the conjugation property~\cite{combettes2011proximal}:\sq
\begin{equation*}
\label{eq:conj_property}
\prox_{\nu F^\star} \bs \zeta = \bs \zeta - \nu \prox_{\frac{1}{\nu}F} \tfrac{1}{\nu} \bs \zeta.\sq
\end{equation*}
The proximal operator of $F_1$ is given by the soft thresholding operator~\cite{combettes2011proximal}\sq
$$ 
\prox_{\frac{1}{\nu} F_1} \bs \zeta = \text{soft}_{[-\frac{1}{\nu},\frac{1}{\nu}]}(\bs \zeta) = \text{sign}(\bs \zeta) (|\bs \zeta| - \tfrac{1}{\nu})_+.\sq 
$$
 
The proximal operators of $F_2$, $F_3$ and $H$ are given by the projection onto the convex sets $\cl C_1$, $\cl C_2$ and $\Omega$, respectively:
\begin{align*}
\prox_{\frac{1}{\nu} F_2} \bs \zeta&= \bs \zeta \min \big(1, \tfrac{\varepsilon_{\rm n}}{ \| \bs \zeta \|_2 }\big),\\  
\prox_{\frac{1}{\nu} F_3} \bs \zeta&= \bs q + \text{soft}_{[-\lambda_{\rm w},\lambda_{\rm w}]}(\bs \zeta - \bs q),\\
\prox_{\mu H} \bs \zeta &= {\rm diag}(0,1,\,\cdots,1)\,\bs \zeta,\sq
\end{align*}
with $\lambda_{\rm w} = 0$ if $\| \bs \zeta - \bs q \|_1 \leq \varepsilon_{\rm w}$ and, otherwise, $\lambda_{\rm w}$ is found by solving $\sum_{i=1}^{2N} \max \{0, |\zeta_i - q_i| - \lambda_{\rm w}\} = \varepsilon_{\rm w}$.

In order to guarantee the convergence of the algorithm, \ie to ensure that $\bs w^{(k)}$ converges to the solution of \eqref{eq:tvl1_2} when $k$ increases, we need to set $\mu$ and $\nu$ such that $\mu \nu \bbb \bs K \bbb^2<1$ \cite{ChambollePock_2011}. The induced norm of the operator ($\bbb \bs K \bbb$) is estimated using the standard power iteration algorithm \cite{sidky2012convex}.

\section{Results}
\label{sec:results}

In this section, we validate our convex approach by studying the quality of the unwrapped phase with respect to the amount of ``wraps'' in the measurements, the noise level and the presence of non-Itoh discontinuities in the original phase image. Results are fairly compared with a post-denoised phase unwrapping obtained by the conjunction of the PUMA algorithm \cite{Bioucas2007} with an optimal soft thresholding denoising \cite{donoho1995} using the same wavelet basis as in \eqref{eq:tvl1_1}. Hereafter, the solutions of our convex approach and of the post-denoised PUMA are denoted as $\bs x_{\rm c}$ and $\bs x_{\rm dp}$, respectively.

Two kinds of discrete phase images are selected in our experiments. They are defined on a $256\!\times\!256$ pixel grid ($N=256^2$). In the first image the phase is simulated by a 2-D Gaussian function of height $0.9\pi$, and standard deviations of 40 pixels horizontally and 25 pixels vertically. In the second image, the phase is simulated by a truncated version of the 2-D Gaussian, where the image is masked by a side triangle. By truncating the Gaussian image, we are able to simulate phase discontinuities.
%\begin{figure}[h]
%  \centering
%  	\includegraphics[width=4cm]{Fig1a}
%	\hspace{0.3cm}
%	\includegraphics[width=4cm]{Fig1b}\sq
%\caption{Example phase images. Left: Gaussian. Right: Truncated Gaussian.\sq}
%\label{fig:Fig1}
%\end{figure}

\begin{figure*}[t]
	\begin{center}
	\includegraphics[width=4.3cm]{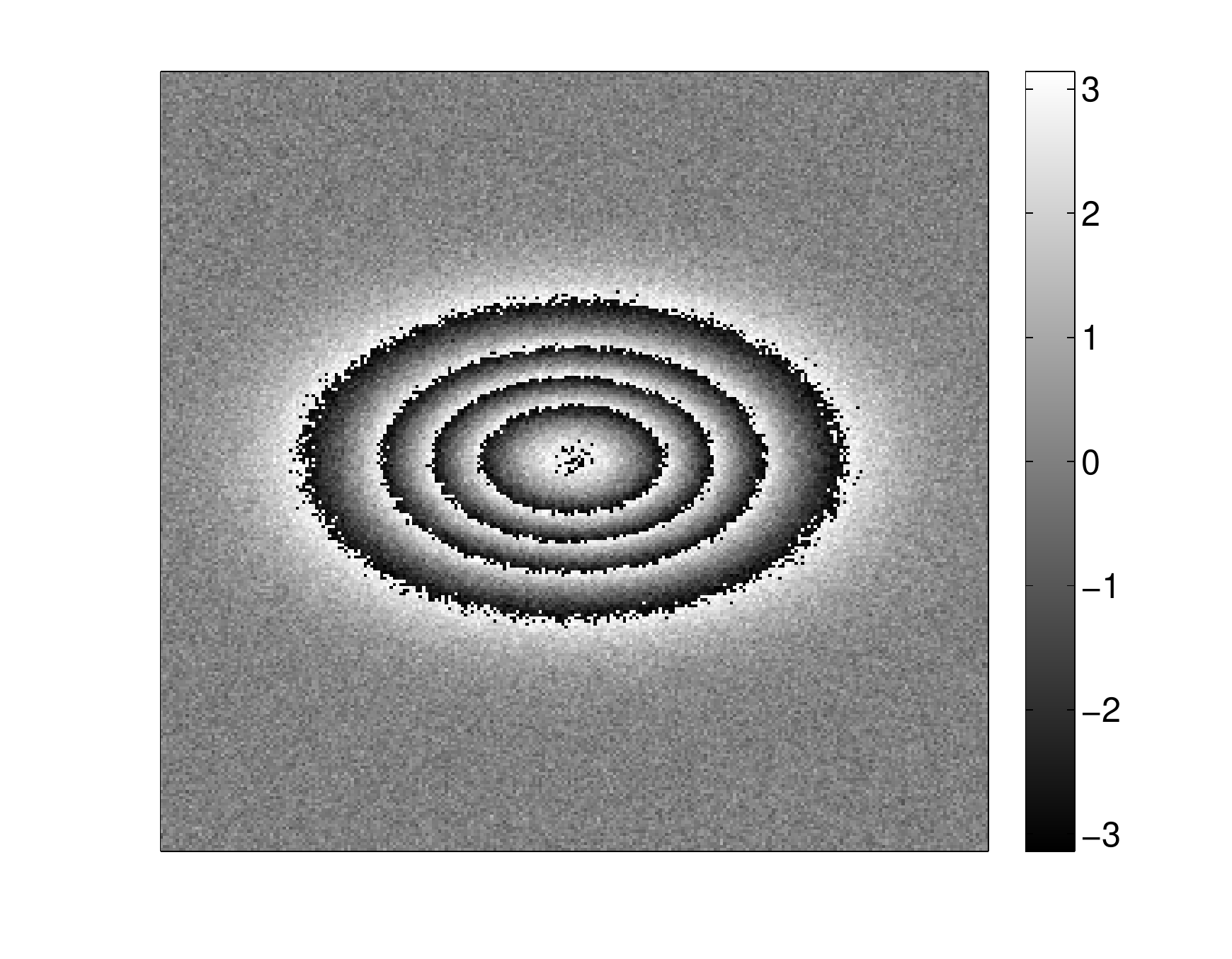}
	\includegraphics[width=4.3cm]{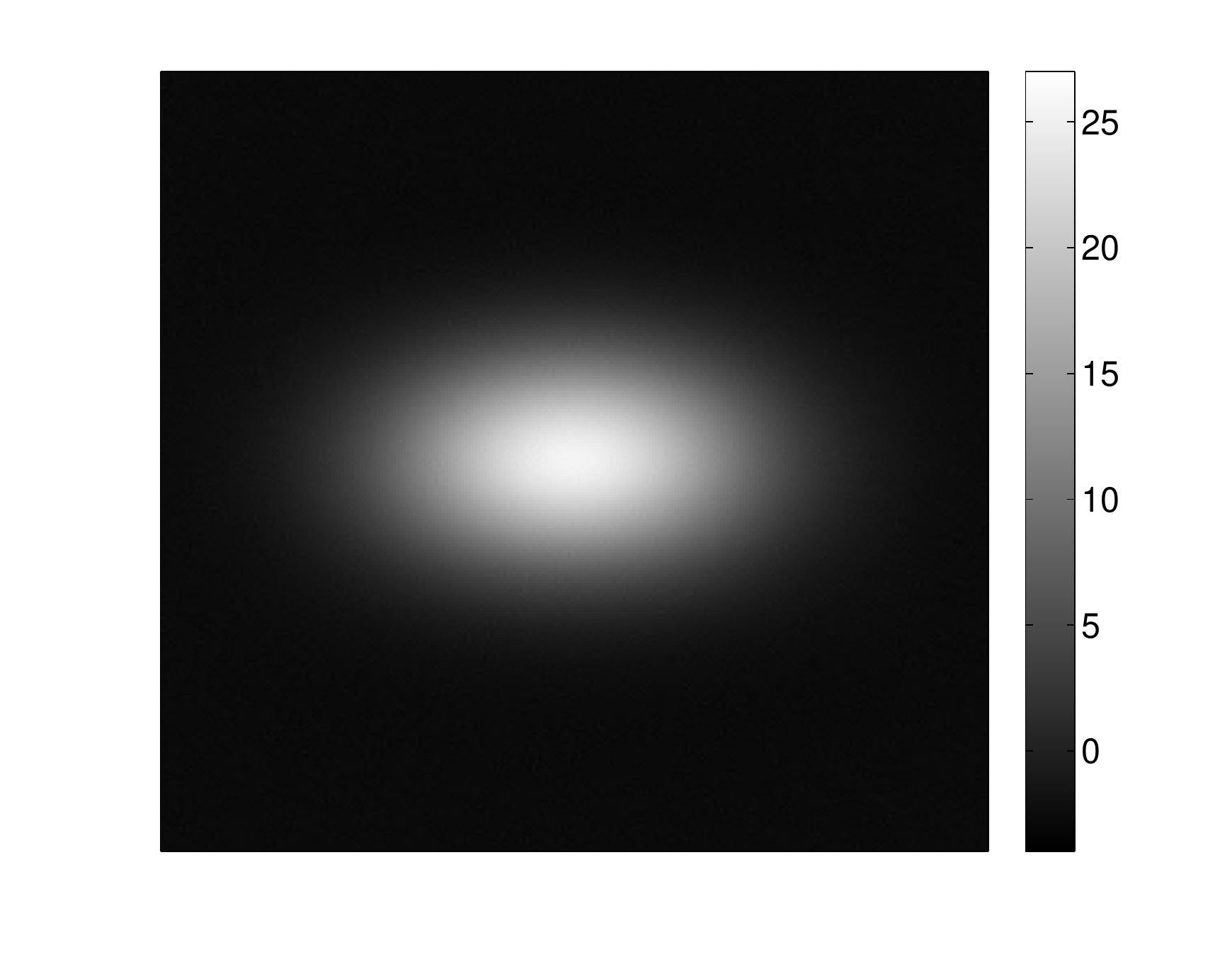}
	\includegraphics[width=4.3cm]{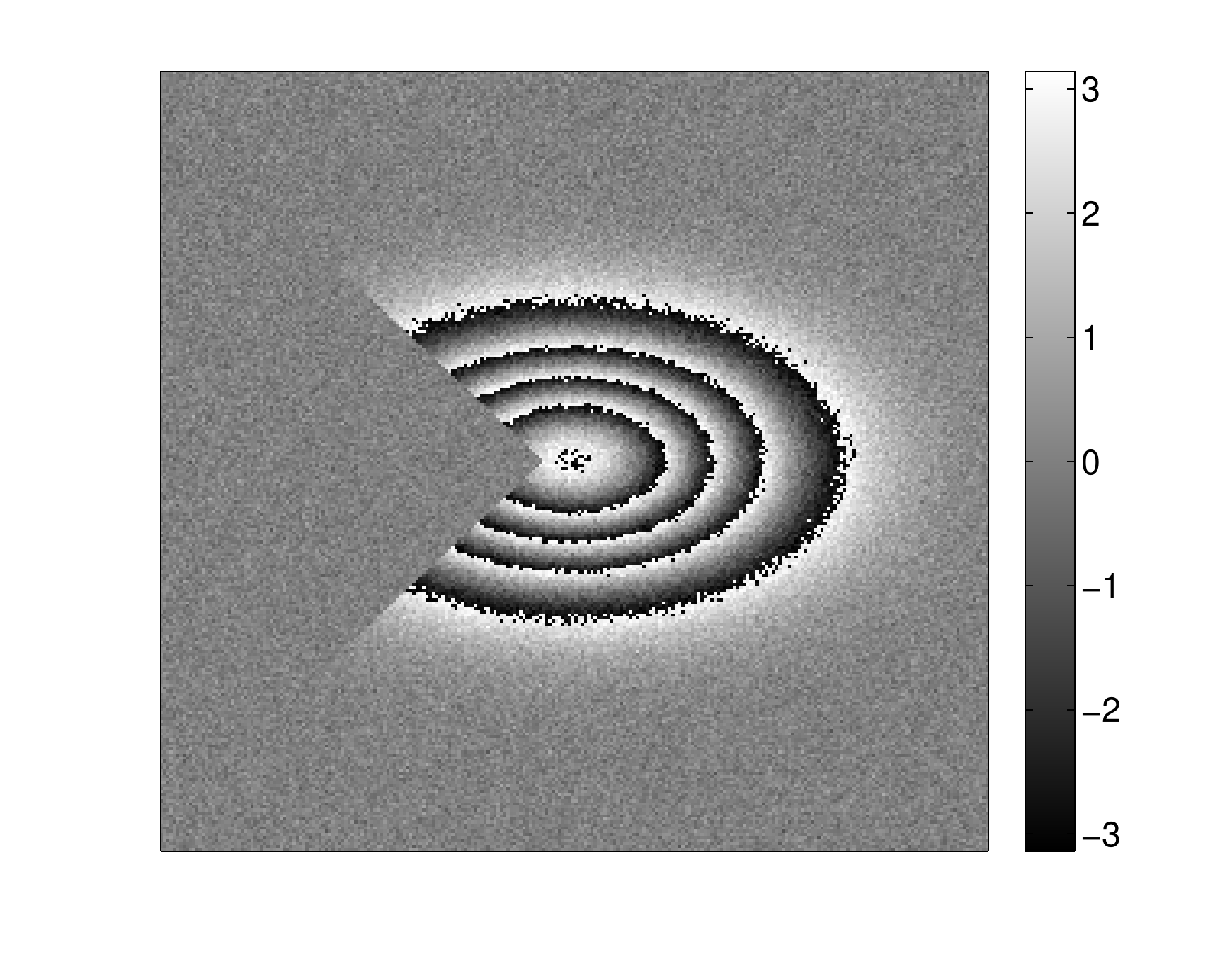}
	\includegraphics[width=4.3cm]{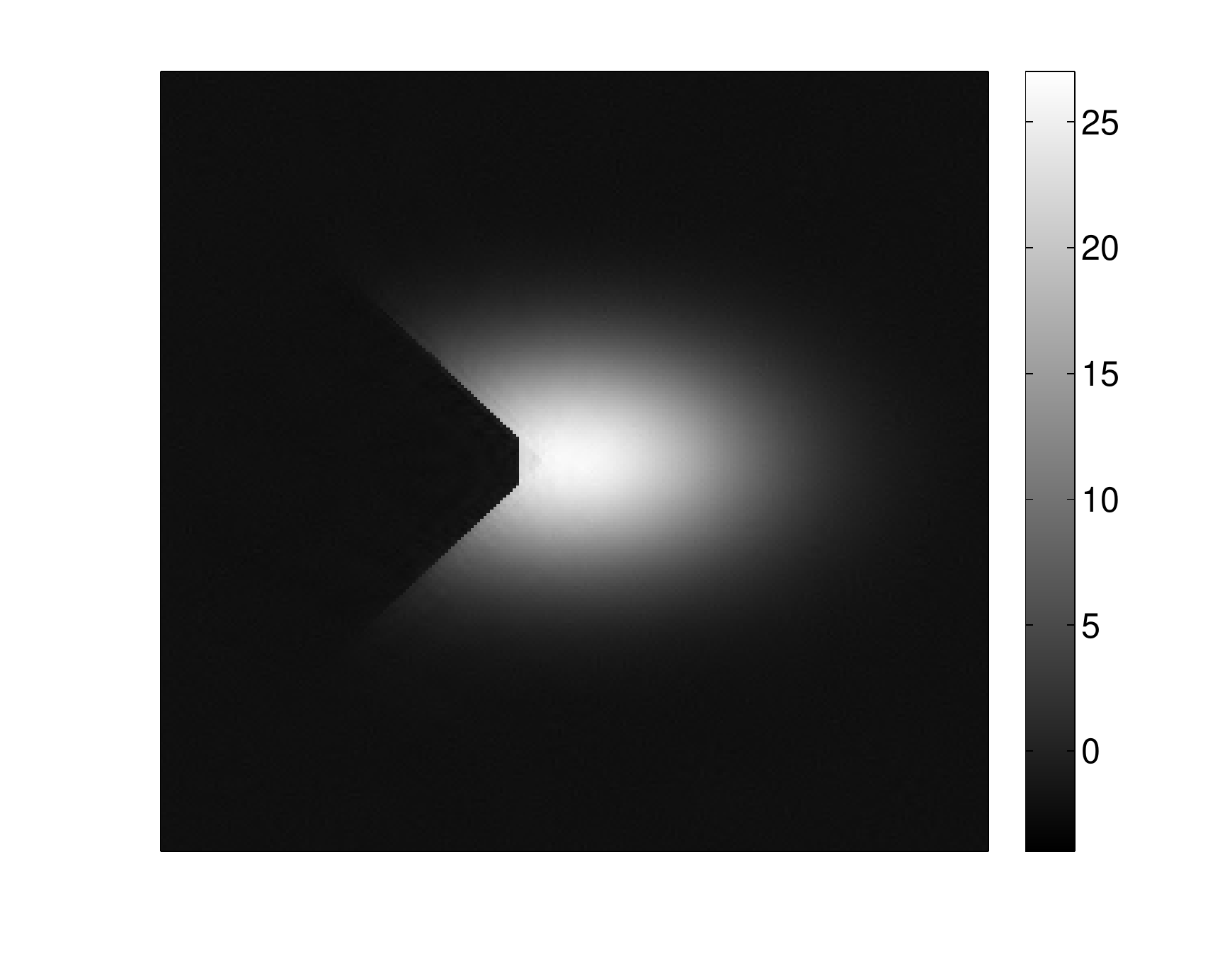}\\
	\end{center}
	\vspace{-0.6cm} \footnotesize {\hspace{1.95cm} (a) \hspace{3.95cm} (b) \hspace{3.95cm} (c) \hspace{3.95cm} (d)}
\caption{Reconstruction results. (a) Noisy measurement of 2-D Gaussian for $\rho = 10$, ISNR = 25dB. (b) Reconstruction of image in (a) using our convex approach (RSNR = 35.1dB). (c) Noisy measurement of 2-D Truncated Gaussian for $\rho = 10$, ISNR = 25dB. (d) Reconstruction of image in (c) using our convex approach (RSNR = 16.7dB).\sq}
\label{fig:Fig1}
\end{figure*}

The robustness of our method is tested against three different noise levels, each characterized by a different ``Input SNR'', $\textrm{ISNR} = 20 \log_{10} \| \bs x \|_2 / \| \bs n \|_2$, namely, 10 dB, 25 dB and $\infty$ dB (no noise). Since there is no reliable estimation of the initial signal mean in any phase unwrapping method, all reconstruction qualities are measured with centralized reconstructions and ground truths, \ie by mean subtraction. After such procedure, the quality of a given reconstruction $\tilde{\bs x} \in \{ \bs x_{\rm c}, \bs x_{\rm dp} \}$ is measured with the ``Reconstruction SNR'': $\textrm {RSNR}\ =\ 20\log_{10} \|\bs x\|_2/\|\bs x - \tilde{\bs x}\|_2$.

All algorithms were implemented in Matlab and executed on a 3.2 GHz CPU, running a 64 Bit Linux system.

For the behavior of the algorithm with respect to the amount of ``wraps'', we analyze the Gaussian phase image and we vary its intensity by multiplying the image by a factor $\rho\in [1, 20]$, $\rho = 1$ providing no phase wraps since $\|\bs x\|_\infty < \pi$. We noticed that for the noiseless scenario (ISNR = $\infty$ dB), the reconstruction quality is not affected by the value of $\rho$, since the Itoh condition is always satisfied. However, since there is no noise, PUMA outperforms the proposed method for all $\rho$. Table \ref{tab:rho_Comparison} presents a comparison for the two noisy scenarios, \ie ISNR = 25dB and ISNR = 10dB. The RSNR is presented for an average of 5 trials. 

\begin{table}[h]
	\centering
	\begin{tabular}{ c | c | c | c | c |}
		\cline{2-5}
		& \multicolumn{4}{|c|}{$\rm RSNR [\rm dB]$}\\
		\cline{2-5}
		 & \multicolumn{2}{|c|}{$\rm ISNR = 25\,{\rm dB}$} & \multicolumn{2}{|c|}{$\rm ISNR = 10\,{\rm dB}$}\\
		\cline{1-5}
		 \multicolumn{1}{|c|}{$\rho$}  & C & DP & C & DP \\ \hline
		\multicolumn{1}{|c|}{1} & 34.18 & 34.11 & 29.62 & 19.14 \\ \hline
		\multicolumn{1}{|c|}{5}  & 42.27 & 34.13 & 22.17 & 20.69 \\ \hline
		\multicolumn{1}{|c|}{10} & 35.08 & 34.13 & 5.04 & 5.05 \\ \hline
		\multicolumn{1}{|c|}{20} & 35.18 & 34.34 & fail & fail  \\ \hline
	\end{tabular}
	\caption{Comparison of the different RSNR obtained using denoised-PUMA (DP) and our convex approach (C) for different values of $\rho$ on the Gaussian phase image.}
  \label{tab:rho_Comparison}
\end{table}

We can notice that the convex approach outperforms the denoised-PUMA (DP) for the scenarios where the Itoh condition is satisfied. However, for those cases where this condition is affected by the noise corrupting the phase, DP provides similar results. About the numerical complexity, for the first noise scenario (ISNR = 25dB), the convex algorithm convergence is reached for an average of 10000 iterations and it takes approximately 11 minutes; while for the second noise scenario (ISNR = 10dB), the convergence is reached for an average of 15000 iterations and takes approximately 16 minutes. 

Fig. \ref{fig:Fig1} depicts the resulting images for the Gaussian and the Truncated Gaussian phases\footnote{Remark that the Ground Truth images are not shown due to lack of space but they are visually very close to (b) and (d).}. Results are shown for $\rho = 10$ and ISNR = 25dB. For the Gaussian phase image, the convex approach provides a good reconstruction quality with RSNR = 35.1dB. We can note that for the Truncated Gaussian the reconstruction quality decreases with RSNR = 16.7dB, because the algorithm is not able to completely recover the phase due to the high discontinuity at the peak of the triangle.

\section{Conclusion}
\label{sec:conclusion}

We propose a general convex optimization approach for robust phase unwrapping. In contrast to state-of-the-art techniques, the proposed approach aims at simultaneously unwrap and denoise the phase image. The proposed approach is shown to outperform the post-denoised PUMA for those scenarios where the noisy phase is smooth enough to satisfy the Itoh condition. However, when such condition is violated due to the noise level or discontinuities in the phase image, the algorithm is not capable of recovering the phase with high quality and it presents the same quality as the denoised PUMA. In future works, we could envisage to remove  from the reconstruction problem the few pixels where the noisy phase is not smooth enough. However the question remains in how to obtain a good estimation on the position of those discontinuities since it depends on the phase to recover.

\end{document}